\documentclass[12pt]{article}

\usepackage{amsthm,amsmath,amssymb}

\usepackage{graphicx}
 \usepackage{epsfig}
  \usepackage{latexsym, epsfig, psfrag,eepic,colordvi,bm}
  \usepackage{graphics,color}
  \usepackage{graphics}
  \usepackage{graphicx}
  \usepackage{epsfig}
  \usepackage[all]{xy}
  
  \usepackage{amssymb}
  \usepackage{amsthm,amsmath}
  
  \usepackage{latexsym, epsfig, psfrag,eepic,colordvi,bm}
  \usepackage{graphics,color}
  \usepackage{amsmath,amsfonts,amssymb,amscd}
  \usepackage[all]{xy}
  \usepackage{epstopdf}
  
  \usepackage{amsthm,amsmath,amssymb}
  
  \usepackage{graphicx}
  \usepackage{amsthm,amsmath,amssymb}
  
  \usepackage{graphicx,epsfig}
\usepackage[colorlinks=true,citecolor=black,linkcolor=black,urlcolor=blue]{hyperref}

\usepackage[all]{xy}

\theoremstyle{plain}
\newtheorem{theorem}{Theorem}[section]

\newtheorem{lemma}{Lemma}[section]
\newtheorem{corollary}{Corollary}[section]
\newtheorem{proposition}{Proposition}[section]

\theoremstyle{definition}
\newtheorem{definition}[theorem]{Definition}
\newtheorem{example}{Example}[section]
\newtheorem{conjecture}[theorem]{Conjecture}

\theoremstyle{remark}
\newtheorem{remark}[theorem]{Remark}

\date{}

\title{\bf On a lower bound for sorting signed permutations by reversals}

\author{Andrei C. Bura$^a$, Ricky X. F. Chen$^b$, Christian M. Reidys$^c$\\
	\small  Biocomplexity Institute and Dept. of Mathematics,\\[-0.8ex]
	\small Virginia Tech, 1015 Life Sciences Circle,\\[-0.8ex]
	\small Blacksburg, VA 24061, USA\\
	\small\tt $^a$anbur12@vbi.vt.edu, $^b$chen.ricky1982@gmail.com, $^c$duck@santafe.edu
}

\begin{document}

\maketitle

\begin{abstract}
	
	Computing the reversal distances of signed permutations is an important topic in Bioinformatics. Recently, a new lower bound for the reversal distance was obtained via the plane permutation framework. This lower bound appears different from the existing lower bound obtained by Bafna and Pevzner through breakpoint graphs. 
	In this paper, we prove that the two lower bounds are equal.
	Moreover, we confirm a related conjecture on skew-symmetric plane permutations, which can be restated as follows:
	let $p=(0,-1,-2,\ldots -n,n,n-1,\ldots 1)$ and let 
	$$
	\tilde{s}=(0,a_1,a_2,\ldots a_n,-a_n,-a_{n-1},\ldots -a_1)
	$$ 
	be any long cycle on the set $\{-n,-n+1,\ldots 0,1,\ldots n\}$. Then, $n$ and $a_n$ are always in the same cycle of the product $p\tilde{s}$. Furthermore, we show the new lower bound 
 via plane permutations can be interpreted as the topological genera of orientable surfaces associated to signed permutations.


\end{abstract}

\section{Introduction}

In Bio-informatics, comparative study of genome sequences is very important to understand evolution.
In particular, the problem of determining the minimum number of certain operations required to transform
one of two given genome sequences into the other, is extensively studied.
Combinatorially, this problem can be formulated as sorting a given permutation (or sequence)
to the identity permutation, by certain operations, in a minimum number of steps, called the distance between the permutation to be sorted and the identity, w.r.t. the operations chosen.
The main operations studied are the transposition~\cite{bul,pev1,chri1,eli}, the block-interchange~\cite{chri1,chri2} and
the reversal~\cite{pev2,pev4,pev3},
 or some combination of them~\cite{baoh, dias}.

In Chen and Reidys~\cite{chr-1}, a framework based on plane permutations was proposed in order to study
these distance problems. As results, they give a general formulation for lower bounds of the transposition and block-interchange
distance from which the existing lower bounds, obtained by Bafna and Pevzner~\cite{pev1}, and Christie~\cite{chri2},
can easily be derived. As to the reversal distance problem of signed permutations,
the authors translate it into a block-interchange distance problem, by which they obtain a new lower bound. They observe that the new lower bound gives the exact reversal
distance for most of the signed permutations.

This paper is mainly concerned with the reversal distances for signed permutations.
A lower bound for the reversal distance was previously obtained in~\cite{pev2,pev4} by studying breakpoint graphs associated to signed permutations.
In~\cite{hch}, a topological framework for sorting signed permutations by reversals was proposed, where the topological genera of these orientable or non-orientable $\pi$-maps associated to signed permutations give a lower bound.
These two existing lower bounds are equal~\cite{hch}.
However, it was not clear which is better, compared to the lower bound in~\cite{chr-1}.
In this paper, we will present another way of associating orientable fatgraphs to signed permutations whose genera give a lower bound for the reversal distance.
Furthermore, we prove all these three approaches provide the same lower bound.

A brief outline of the paper is as follows: in Section~$2$, we review the approach to study the reversal distance
via plane permutations proposed in~\cite{chr-1}, as well as a conjecture on skew-symmetric plane permutations to be confirmed in this paper.
In Section~$3$, from each signed permutation, we construct a fatgraph and show that the genus of the fatgraph gives a lower bound for the reversal distance of the signed permutation. In Section~$4$, we will show
the equivalence of the lower bound on the reversal distance obtained by Chen and Reidys~\cite{chr-1}  and the lower bound obtained by Bafna and Pevzner~\cite{pev1} through breakpoint graphs as well as the genus bound in Huang and Reidys~\cite{hch}.

\section{Plane permutations and the reversal distance}

\begin{definition}[Plane permutation]
A plane permutation on $[n]=\{1,2,\ldots n\}$ is a pair $\mathfrak{p}=(s,\pi)$ where $s=(s_i)_{i=0}^{n-1}$
is an $n$-cycle and $\pi$ is an arbitrary permutation. The permutation $D_{\mathfrak{p}}=
s\circ \pi^{-1}$ is called the diagonal of $\mathfrak{p}$.
\end{definition}\label{2def1}

Let $s=(s_0,s_1,\cdots,s_{n-1})$.
A plane permutation $(s,\pi)$ can be represented by two aligned rows
\begin{equation}
(s,\pi)=\left(\begin{array}{ccccc}
s_0&s_1&\cdots &s_{n-2}&s_{n-1}\\
\pi(s_0)&\pi(s_1)&\cdots &\pi(s_{n-2}) &\pi(s_{n-1})
\end{array}\right)
\end{equation}
Note that $D_{\mathfrak{p}}$ is determined by the diagonal pairs, cyclically, in the two-line
representation, i.e., $D_{\mathfrak{p}}(\pi(s_{i-1}))=s_i$ for $0<i< n$, and
$D_{\mathfrak{p}}(\pi(s_{n-1}))=s_0$.
In the following, by ``the cycles of $\mathfrak{p}=(s,\pi)$'' we mean the cycles of $\pi$.

Given a plane permutation $\mathfrak{p}=(s,\pi)$ on $[n]$ and a sequence $h=(i,j,k,l)$, such that $i\leq j<k \leq l$
and $\{i,j,k,l\}\subset [n-1]$, let
$$
s^h=(s_0,s_1,\dots,s_{i-1},\underline{s_k,\dots,s_l},s_{j+1},\dots,s_{k-1},
\underline{s_i,\dots,s_j},s_{l+1},\dots).
$$
In other words $s^h$ is the $n$-cycle obtained by transposing the blocks $[s_i,s_j]$ and $[s_k,s_l]$ in $s$.
Let furthermore
$$
\pi^h=D_{\mathfrak{p}}^{-1}\circ s^h.
$$
This means, the derived plane permutation $\mathfrak{p}^h=(s^h,\pi^h)$ can be represented as

\begin{eqnarray*}
\left(
\vcenter{\xymatrix@C=0pc@R=1pc{
\cdots s_{i-1}  & s_k\ar@{--}[dl] &\cdots & s_l\ar@{--}[dl] & s_{j+1} &\hspace{-0.5ex}\cdots\hspace{-0.5ex} &
s_{k-1} & s_i\ar@{--}[dl] &\cdots & s_{j}\ar@{--}[dl]& s_{l+1}  \cdots\\
\cdots \pi(s_{k-1}) & \pi(s_k) & \cdots & \pi(s_j) & \pi(s_{j+1}) &\hspace{-0.5ex}\cdots\hspace{-0.5ex} & \pi(s_{i-1}) & \pi(s_i) &\cdots & \pi(s_l)& \pi(s_{l+1}) \cdots
}}
\right).
\end{eqnarray*}
Note that the bottom row of the two-row representation of $(s^h,\pi^h)$
is obtained by transposing the blocks $[\pi(s_{i-1}),\pi(s_{j-1})]$ and
$[\pi(s_{k-1}),\pi(s_{l-1})]$ of the bottom row of $(s,\pi)$.
In this way, each $h$ not only determines a block-interchange on the $s$-sequence, but also it induces a block-interchange on the plane permutation (in two-row representation).

\begin{definition}
A signed permutation on $[n]$ is a pair $(a,w)$ where $a=a_1a_2\cdots a_n$ is
a sequence on $[n]$ while $w=w_1w_2\cdots w_n$ is a word of length $n$ on the
alphabet set $\{+,-\}$.
\end{definition}
Usually, a signed permutation is represented by a single sequence
$a_w=a_{w,1}a_{w,2}\cdots a_{w,n}$ where $a_{w,k}=w_ka_k$, i.e.,
each $a_k$ carries a sign determined by $w_k$.

Given a signed permutation $a=a_1a_2\cdots a_{i-1}a_i a_{i+1}\cdots a_{j-1}
a_j a_{j+1}\cdots a_n$ on $[n]$, a reversal $\varrho_{i,j}$ acting on
$a$ will change $a$ into
$$
a'=\varrho_{i,j}\diamond a =a_1a_2\cdots a_{i-1}(-a_j)(-a_{j-1})\cdots
(-a_{i+1})(-a_i)a_{j+1}\cdots a_n.
$$
The reversal distance $d_r(a)$ of a signed permutation $a$ on $[n]$ is the minimum
number of reversals needed to sort $a$ into $e_n=12\cdots n$.\

\begin{example}\label{2exam1}
The signed permutation $-5,1,-3,2,4$ needs at least $4$ steps to be sorted as illustrated below: 
\begin{eqnarray*}
-5 \quad +1 \quad \underline{-3 \quad +2} \quad +4\\
-5 \quad +1 \quad \underline{-2} \quad +3 \quad +4\\
\underline{-5 \quad +1 \quad {+2} \quad +3 \quad +4}\\
\underline{-4 \quad -3 \quad -2 \quad -1} \quad +5\\
+1 \quad +2 \quad +3 \quad +4 \quad +5
	\end{eqnarray*}
\end{example}
For a given signed permutation $a$, we associate to it the sequence $s=s(a)$ as follows
$$
s=s_0s_1s_2\cdots s_{2n}=0a_1a_2\cdots a_n (-a_n)(-a_{n-1})\cdots (-a_2)(-a_1).
$$
In other words $s_0=0$ and $s_{k}=-s_{2n+1-k}$ for $1\leq k\leq 2n$. Furthermore, such sequences
will be referred to as skew-symmetric sequences since we have $s_{k}=-s_{2n+1-k}$. A sequence
$s$ is called negative if there exists $s_i<0$ for some $1\le i\le n$.
Now we observe that the reversal distance of $a$ is equal to the block-interchange distance of $s$
into
$$
e_n^{\natural}=012\cdots n (-n)(-n+1)\cdots (-2) (-1),
$$
where only certain block-interchanges are allowed. In other words, only
block-interchanges $h=(i,j,2n+1-j,2n+1-i)$ where $1\leq i\leq j\leq n$  are allowed. Hereafter, we will denote
these particular block-interchanges on $s$ as the reversals $\varrho_{i,j}$.

Let
\begin{eqnarray*}
\tilde{s}&=&(s)=(0,a_1,a_2,\ldots,a_{n-1},a_n,-a_n,-a_{n-1},\ldots,-a_2,-a_1),\\
p&=&(0,-1,-2,\ldots,-n+1,-n,n,n-1,\ldots,2,1).
\end{eqnarray*}
A plane permutation of the form $(\tilde{s},\pi)$ will be called skew-symmetric.
Let $C(\pi)$ denote the number of cycles in the permutation $\pi$.
Based on these definitions and notations, the new lower bound for the reversal distance of $a$ obtained in~\cite{chr-1}
is

\begin{theorem}[Chen and Reidys~\cite{chr-1}]\label{7thm1}
\begin{equation}\label{5eq1}
d_r(a)\geq \frac{2n+1-C(p\tilde{s})}{2}.
\end{equation}
\end{theorem}

The proof of the theorem relies on the observation that
the increment of the number of cycles from $\mathfrak{p}$ to $\mathfrak{p}^h$,
for a given reversal $h$, is at most $2$.
Now start with $\mathfrak{p}=(\tilde{s},\pi_1)$ such that $D_{\mathfrak{p}}=p^{-1}$,
and let us apply a sequence of reversals so that we arrive at $((e_n^{\natural}),\pi_2)$. Note that
the starting plane permutation
has $C(p\tilde{s})$ cycles, while
the terminating plane permutation has $C(p(e_n^{\natural}))=2n+1$ cycles.
Hence, at least $\frac{2n+1-C(p\tilde{s})}{2}$ reversals are needed for sorting ${s}$ into $e_n^{\natural}$.

The effectiveness of the lower bound, when there exists a reversal increasing the number of
cycles by exactly $2$ (a $2$-reversal for short), was analyzed in~\cite{chr-1}.
It was proved that there exists a $2$-reversal for all skew-symmetric plane permutations
but possibly two minor classes. One class is that of nonnegative permutations.
For the other class, a $2$-reversal exists if the following conjecture is true.


\begin{conjecture}\label{7conj2}\cite{chr-1}
Let $\mathfrak{p}=(\tilde{s},\pi)$ be a skew-symmetric plane permutation on $[n]^{\pm}=\{-n,-n+1,\ldots 0,1,\ldots n\}$,
where $D_{\mathfrak{p}}=p^{-1}$.
Then, $n$ and $s_n$ are in the same cycle of $\pi$.
\end{conjecture}

\section{Associating fatgraphs to signed permutations}

Let $G$ be a graph, with loops and multiple edges allowed.
We call the ends of edges in $G$ half-edges (or darts).
An orientable fatgraph, also called {\it (combinatorial) map}, is a graph
$G$ with a specified cyclic order of half-edges around (i.e., incident to) each vertex of $G$.
An orientable fatgraph can be viewed as a cell-complex, i.e., the underlying graph represents the $1$-skeleton while the (counterclockwise) cyclic order of half-edges indicate how to glue $2$-cells on top of that to obtain a cell-complex.
Accordingly, we can define the topological genus of the fatgraph to be the genus of the cell-complex. Hereafter, a fatgraph is alway orientable unless explicitly stated otherwise.

A fatgraph having $n$ edges can be represented as a triple of permutations $(\alpha,\beta,\gamma)$ on
$[2n]$ where $\alpha$ is a fixed-point free involution and $\gamma=\alpha\beta$.
This can be seen as follows: we 
label the half-edges of the fatgraph using the labels from the set $[2n]$ so that each label appears exactly
once. This induces two permutations $\alpha$ and $\beta$, where $\alpha$ is a fixed point
free involution, whose cycles consist of the labels of the two half-edgs of the same
(untwisted) edge and $\beta$-cycles represent the counterclockwise cyclic arrangement of
all half-edges incident to the same vertex. $\gamma=\alpha\beta$-cycles are called the
{\it faces}. The topological genus of a fatgraph $(\alpha,\beta,\gamma)$ satisfies
\begin{align}
C(\beta)-C(\alpha)+C(\gamma)=2-2g.
\end{align}
Conversely, from such a triple of permutations $(\alpha,\beta,\gamma)$ on
$[2n]$, a fatgraph can be constructed.

For the general case that $\alpha$ is not necessarily a fixed-point free involution, the pair $(\alpha,\beta,\gamma)$ satisfying $\gamma=\alpha\beta$ is called a hypermap.

Note a plane permutation $\mathfrak{p}=(s, \pi)$ can be interpreted as a hypermap $(D_{\mathfrak{p}}, \pi, s)$~\cite{chr-1, chr-2} and 
there is a connection between hypermaps and certain maps (called bipartite maps) 
observed by Walsh~\cite{walsh}. 
Motivated by this, we will associate an orientable fatgraph for a given signed permutation by transforming the skew-symmetric plane permutation for a signed permutation into a plane permutation whose diagonal is a fixed-point free involution. 
We next show that the topological genus of the fatgraph encoded by the plane permutation gives the lower bound for the reversal distance.

Given a plane permutation $\mathfrak{p}$ on $[n]$, w.l.o.g., we assume
\begin{equation*}
\mathfrak{p}=(s,\pi)=\left(\begin{array}{ccccc}
1&2&\cdots &n-1&n\\
\pi(1)&\pi(2)&\cdots &\pi(n-1) &\pi(n)
\end{array}\right).
\end{equation*}
We obtain another plane permutation $(\hat{s},\hat{\pi})$ on $[n]\cup [\bar{n}]$ having a fixed-point free involution diagonal, where $[\bar{n}]=\{\bar{1},\bar{2},\ldots \bar{n}\}$, by the following construction:
\begin{description}
	\item[step~$1$.] put $\bar{i}$ after $i$ to obtain $\hat{s}$;
	\item[step~$2$.] we determine $\hat{\pi}$ as follows: firstly, when restricted to $[n]$, $\hat{\pi}|_{[n]}=\pi$. Next, to guarantee
	the diagonal to be a fixed-point free involution, the image of $\bar{i}$ (i.e., $\hat{\pi}(\bar{i})$) will be automatically determined. Namely, $\hat{\pi}(\bar{1})=\overline{(\pi^{-1}(2))}$, $\hat{\pi}(\bar{2})=\overline{(\pi^{-1}(3))}$, etc. In general, $\hat{\pi}(\bar{i})=\overline{(\pi^{-1}(i+1))}$ (in the
	sense of modulo $n$). 
\end{description}   

The following is an example to illustrate this construction, where $n=6$,
\begin{equation*}
\mathfrak{p}=\left(\begin{array}{cccccc}
1&2&3&4&5&6\\
4&5&6&1&2&3
\end{array}\right)\Longleftrightarrow
\left(\begin{array}{cccccccccccc}
1&\bar{1}&2&\bar{2}&3&\bar{3}&4&\bar{4}&5&\bar{5}&6&\bar{6}\\
4&\bar{5}&5&\bar{6}&6&\bar{1}&1&\bar{2}&2&\bar{3}&3&\bar{4}
\end{array}\right)=(\hat{s},\hat{\pi}).
\end{equation*}

\begin{lemma}
	When restricted to $[\bar{n}]$,
$\hat{\pi}|_{[\bar{n}]}$ has the same number of cycles as
$D_{\mathfrak{p}}$.	
\end{lemma}
\proof

 Note, identifying $\pi(i)$ and $\bar{i}$, $\hat{\pi}(\bar{i})$ and $i+1$, i.e., two elements on a same diagonal-pair, will preserve cycle structure in the sense that
$$
\hat{\pi}:\bar{i}\rightarrow \hat{\pi}(\bar{i}), \quad D_{\mathfrak{p}}: \pi(i)\rightarrow i+1.
$$
Therefore, $\hat{\pi}|_{[\bar{n}]}$ has the same number of cycles as
$D_{\mathfrak{p}}$. 
\qed

For the example given above, we can check that $D_{\mathfrak{p}}=(153)(264)$ has two cycles, while
$\hat{\pi}|_{[\bar{5}]}=(\bar{1}\bar{5}\bar{3})(\bar{2}\bar{6}\bar{4})$ has also two cycles.

Applying above construction to the associated skew-symmetric plane permutation $\mathfrak{p}=(\tilde{s},\pi)$ to a signed permutation $a$, where $D_{\mathfrak{p}}=p^{-1}$, we can associate 
a fatgraph $\mathcal{F}_a$ to $a$. Let $g(\mathcal{F}_a)$ denote
the genus of the fatgraph $\mathcal{F}_a$.

\begin{theorem}
	The reversal distance $d_r(a)$ for $a$ satisfies
	\begin{align}\label{top-bound}
	d_r(a) \geq g(\mathcal{F}_a).
	\end{align} 
\end{theorem}
\proof Let $\mathfrak{p}=(s, \pi)$ be the corresponding skew-symmetric plane permutation of $a$. Then, in $\mathcal{F}_a$, there are $C(\pi)+C(D_{\mathfrak{p}})=C(ps)+1$ 
vertices, $2n+1$ edges and one face. Thus, 
$$
2-2g(\mathcal{F}_a)=(C(ps)+1)-(2n+1)+1
$$
so that
$$
g(\mathcal{F}_a)=\frac{2n+1-C(ps)}{2}.
$$
Comparing with Theorem~\ref{7thm1} completes the proof. \qed

\begin{example}
Consider the signed permutation in Example~\ref{2exam1}.
Its associated skew-symmetric plane permutation is
$$
\left(\begin{array}{ccccccccccc}
	0& -5& 1& -3& 2& 4& -4& -2& 3& -1& 5\\
	5& 0& -4& 1& 3 & -5& -3& 2& -2& 4 & -1
\end{array}\right).
$$
Applying the transformation above, we obtain
$$
\left(\begin{array}{cccccccccccccccccccccc}
0& \bar{0} & -5& \overline{-5} & 1& \cdots & \bar{2} & 4&\bar{4} & -4&\overline{-4} & -2&\overline{-2} & \cdots & -1&\overline{-1} & 5 &\bar{5} \\
5& \bar{4}& 0& \overline{-3}& -4 & \cdots& \overline{-1} & -5 & \bar{1}& -3& \bar{3}& 2& \bar{2} & \cdots & 4& \bar{0}& -1 & \overline{-5}
\end{array}\right).
$$
and its associated fatgraph as illustrated in Figure~\ref{2fig1}. It can be seen that the fatgraph has $4$ vertices, $11$ edges (ribbons) and one face. So, the genus of the fatgraph is $4$ by Euler's characteristic formula.
\begin{figure}[!htb]
	\centering
	\includegraphics[scale=.45]{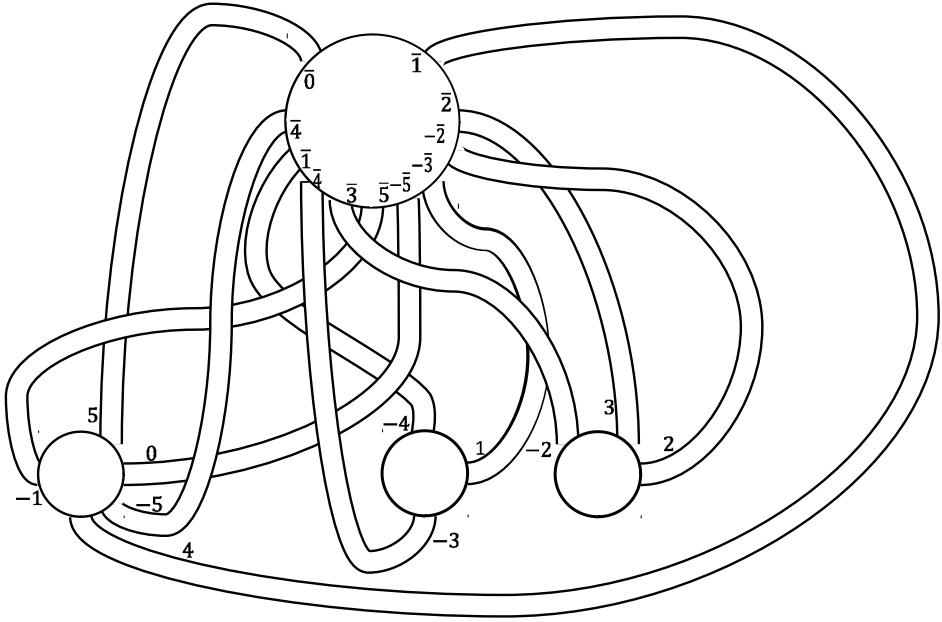}
	\caption{The fatgraph associated to $-5,1,-3,2,4$.}\label{2fig1}
\end{figure}
\end{example}

Compared to the fatgraphs associated to signed permutations in~\cite{hch}, the associated fatgraphs there are non-orientable in most cases while the associated fatgraphs in this paper are always orientable. The relation between the two fatgraphs (via these two different approaches) associated to the same signed permutation remains unclear and would be an interesting 
problem to investigate. 

\emph{Remark.} For the lower bounds for the transposition distances and the block-interchange distances
of (unsigned) permutations, we can have analogous topological interpretations by 
transforming the associated plane permutations~\cite{chr-1} into fatgraphs in the same way as we did
here for the reversal distances of signed permutations.

\section{Comparison with the Bafna-Pevzner lower bound}

In~\cite{hch}, it was shown that the topological-genus lower bound there and the lower bound given by breakpoint graphs in~\cite{pev2,pev4} are equal.
In this section, we will show that our new topological-genus lower bound equals the lower bound given by breakpoint graphs in~\cite{pev2,pev4} as well.

The breakpoint graph for a given signed permutation $a=a_1a_2\cdots a_n$ on $[n]$ can be obtained as follows:
replace $a_i$ with $(-a_i) a_i$, and add $0$ at the beginning of the obtained sequence while
adding $-(n+1)$ at the end. In this way we obtain a sequence $b=b_0b_1b_2\cdots b_{2n}b_{2n+1}$ on $[n]^{\pm} \cup \{-n-1\}$, where $b_0=0, \: b_{2n+1}=-n-1$.
Draw a black edge between $b_{2i}$ and $b_{2i+1}$, as well as a grey edge between $i$ and $-(i+1)$
for $0\leq i \leq n$. The obtained graph is the breakpoint graph, $BG(a)$, of $a$.
Note that each vertex in $BG(a)$ has degree two hence the graph can be decomposed into disjoint cycles.
Denote the number of cycles in $BG(a)$ by $C_{BG}(a)$. Then, the lower bound of the reversal distance via the breakpoint graph is 
\begin{align}\label{bg-r}
	d_r(a) \geq n+1-C_{BG}(a).
\end{align}
In~\cite{chr-1}, an algebraic, equivalent formulation, of the breakpoint graph approach was obtained as follows:
let $\theta_1$,~$\theta_2$ be the two involutions (without fixed points) determined by
the black edges and grey edges in the breakpoint graph, respectively, i.e.
\begin{align*}
	\theta_1 &=(b_0,b_1)(b_2,b_3)\cdots (b_{2n},b_{2n+1}),\\
	\theta_2 &=(0,-1)(1,-2)\cdots (n,-n-1).
\end{align*}
Then, we have
\begin{proposition}[Chen and Reidys~\cite{chr-1}]
	\begin{align}\label{bp-r2}
		d_r(a)\geq \frac{2n+2- C(\theta_1 \theta_2)}{2}.
	\end{align}
\end{proposition}

To show the equivalence between our new topological-genus bound (equivalently, the lower bound in eq.~\eqref{5eq1}) and the lower bound given by the breakpoint graphs, it suffices to show, for any signed permutation $a$, there holds
$$
C(p_r \tilde{s})= C(\theta_1\theta_2)-1.
$$


\begin{definition}
Let $\sigma$ be a permutation on the set $[n]^{\pm}$. We associate to $\sigma$ the matrix $ A_{\sigma}=[a_{ij}]$,
$$
a_{ij}=\left\{
\begin{array}{ll}
1  & \mbox{if } i = \sigma(j) \\
0 & \mbox{if } i \ne \sigma(j)
\end{array}
\right.$$
where $i,j\in [n]^{\pm}$, index rows and columns, following the order $-n,-n+1,\ldots -1,0,1,\ldots,n$, respectively. We call $A_{\sigma}$ the permutation matrix associated to $\sigma$ and denote this by $\sigma\sim A_{\sigma}$.
\end{definition}
Recall that $\sigma\sim A_{\sigma}\Leftrightarrow \sigma^{-1}\sim A_{\sigma}^{T}$, that if $\tau\sim A_{\tau}$ then $\sigma\circ\tau=A_{\sigma}A_{\tau}$, and that for $id$, the identity permutation, we have that $id\sim A_{id}=I_{2n+1}$ (i.e. the identity matrix). 


\begin{lemma}
	Let
	\begin{align*}
		p &=(0,-1,-2,\ldots -n,n,n-1,\ldots 1)\sim A_p=P,\\
		\tilde{s}&=(0,a_1,a_2,\ldots a_n,-a_n,-a_{n-1},\ldots -a_1)\sim A_{\tilde{s}}=S.
	\end{align*}
	Let
	$R=[r_{ij}]$ be the $(2n+1)\times (2n+1)$ unitary anti-diagonal matrix, also known as the exchange matrix. Namely,
	$$r_{ij}=\left\{
	\begin{array}{ll}
	1  & \mbox{if } j = 2n-i+2 \\
	0 & \mbox{if } j \ne 2n-i+2
	\end{array}
	\right.$$
	where $i,j\in [2n+1]$ index rows and columns of $R$ respectively.\\
	Then
	\begin{align}
		PS=(PR)(RS),
	\end{align}
\end{lemma}
\proof It suffices to check that $R=R^T$. Then, since $R$ is a permutation matrix we have $R^2=RR^T=I_{2n+1}$.  \qed

 We next show that both permutations $PR$ and $RS$ are involutions with a unique fixed point.

\begin{lemma}
The permutation corresponding to the matrix $PR$ is the involution
$$
p_{invo}=(-n, n-1)(-n+1,n-2)\cdots (-1,0)(n).
$$
\end{lemma}
\proof
We compute
\[PR=\begin{bmatrix}
	0 & 1 & 0 & 0 & \dots  & 0 & 0 \\
	0 & 0 & 1 & 0 & \dots  & 0 & 0 \\
	0 & 0 & 0 & 1 & \dots  & 0 & 0 \\
	\vdots & \vdots & \vdots & \vdots & \ddots & \vdots & \vdots \\
	0 & 0 & 0 & 0 & \dots  & 1 & 0 \\
	0 & 0 & 0 & 0 & \dots  & 0 & 1 \\
	1 & 0 & 0 & 0 & \dots  & 0 & 0 \\
\end{bmatrix}
\begin{bmatrix}
0 & 0 & 0 & 0 & \dots  & 0 & 1 \\
0 & 0 & 0 & 0 & \dots  & 1 & 0 \\
0 & 0 & 0 & 0 & \dots  & 0 & 0 \\
\vdots & \vdots & \vdots & \vdots & \ddots & \vdots & \vdots \\
0 & 0 & 1 & 0 & \dots  & 0 & 0 \\
0 & 1 & 0 & 0 & \dots  & 0 & 0 \\
1 & 0 & 0 & 0 & \dots  & 0 & 0 \\
\end{bmatrix}
=
\begin{bmatrix}
0 & 0 & 0 &  \dots 0 & 1 & 0 \\
0 & 0 & 0 &  \dots 1 & 0 & 0 \\
0 & 0 & 0 &  \dots 0 & 0 & 0 \\
\vdots & \vdots & \vdots & \ddots & \vdots & \vdots\\
0 & 1 & 0 &  \dots 0 & 0 & 0 \\
1 & 0 & 0 &  \dots 0 & 0 & 0 \\
0 & 0 & 0 &  \dots 0 & 0 & 1 \\
\end{bmatrix}
\]
Converting the matrix into a permutation completes the proof. \qed\\

\begin{lemma}
The permutation corresponding to the matrix $RS$ is the involution
$$
s_{invo}= (0,-s_1)(s_1,-s_2)\cdots(s_{n-1},-s_n)(s_n)
$$
\end{lemma}
\proof
It is easy to check that left multiplication by the exchange matrix $R$ reverses the order on the rows of the multiplied matrix $S=[a_{ij}]$, $i,j\in [n]^{\pm}$. 
Since $a_{i,j}=1\Longleftrightarrow i=s(j)$, when multiplied by the matrix $R$, the row indexed by $i$, in the matrix $S$, is sent to the row indexed by $-i$. This is due to the symmetry of $[n]^{\pm}$. However, this in turn means that in the permutation  $s_{invo}\sim RS$, $j\xrightarrow{s_{invo}}-i,\:\:\forall j\in [n]^{\pm}\:\text{ such that }\:s(j)=i$. But now by virtue of the structure of the cycle $s=(0,s_1,s_2,\ldots s_n,-s_n,-s_{n-1},\ldots -s_1)$, we have 
$$
s_{invo}= (0,-s_1)(s_1,-s_2)\cdots(s_{n-1},-s_n)(s_n)
$$ 
and the lemma follows.
\qed\\\\
We now know that both $PR$ and $RS$ are matrices corresponding to involutions, each with a unique fixed point $n$ and $s_n$,
respectively. We proceed to prove some results regarding the product of two involutions.

 \begin{lemma}\label{noodd}
 	Let $\sigma_1,\: \sigma_2$ be two fixed-point free involutions on a set $T$, with $|T|= 2n$ for $n\geq 1$.
 	Then, there does not exist $k\geq 1$ such that $\sigma_1(\sigma_2\sigma_1)^k(x)=x$.
 \end{lemma}
 \proof 
 Assume by contradiction that there exists a $k\geq 1$ such that $\sigma_1(\sigma_2\sigma_1)^k(x)=x$. Then, since $\sigma_1$ is an involution, we have $\sigma_1(x)=(\sigma_2\sigma_1)^k(x)$.

 Now if $k=1$, we have $\sigma_1(x)=\sigma_2\sigma_1(x)$, implying that $\sigma_1(x)$ is a fixed point of $\sigma_2$, which is a contradiction.

 If $k>1$, 
 since $\sigma_2$ is an involution, we have
 $$
 \sigma_2\sigma_1(x)=(\sigma_1\sigma_2)^{k-1}(\sigma_1(x))=\sigma_1(\sigma_2\sigma_1)^{k-1}(x).
 $$
 Now, if $k-1=1$, we have that $\sigma_2\sigma_1(x)$ has to be a fixed point of $\sigma_1$, which is again a contradiction. Otherwise, we can set $y=\sigma_2\sigma_1(x)$ and obtain
 $y=\sigma_1(\sigma_2\sigma_1)^{k-1}(y)$, and iterate the previous argument.
 In this way, we eventually obtain a fixed point, either for $\sigma_1$ or $\sigma_2$, which contradicts our assumption, hence the lemma follows. \qed\\

 \begin{lemma}\label{keylem}
 	Let $\sigma_1,\:\sigma_2$ be two involutions on a set $T$
 	such that each of them has a unique fixed point, $a$ and $b$ respectively.
 	Then $\sigma_2 \sigma_1$ has a cycle which contains both $a$ and $b$.
 \end{lemma}
 \proof
 There is nothing to prove if $a=b$, so we will assume $a\neq b$ in the following.
 Let now
 \begin{align*}
 	\sigma_1 &=(p_1, p_2)(p_3,p_4) \cdots (p_{2t-1}, p_{2t})(a),\\
 	\sigma_2 &=(q_1, q_2)(q_3,q_4) \cdots (q_{2t-1}, q_{2t})(b),\\
 	\sigma'_1 &=(p_1, p_2)(p_3,p_4) \cdots (p_{2t-1}, p_{2t})(a,x),\\
 	\sigma'_2 &=(q_1, q_2)(q_3,q_4) \cdots (q_{2t-1}, q_{2t})(b,x),
 \end{align*}
 where $\sigma'_1$ and $\: \sigma'_2$ are involutions on $T\cup \{x\}$, where $x\notin T$. We now compare the following two iterations:
 \begin{align*}
 	& a\rightarrow \sigma_1(a)\rightarrow \sigma_2\sigma_1(a) \rightarrow \sigma_1\sigma_2\sigma_1(a)\rightarrow (\sigma_2\sigma_1)^2(a)\cdots (\sigma_2\sigma_1)^{k_1}(a)=a,\\
 	& x\rightarrow \sigma'_1(x)\rightarrow \sigma'_2\sigma'_1(x) \rightarrow \sigma'_1\sigma'_2\sigma'_1(x)\rightarrow (\sigma'_2\sigma'_1)^2(x)\cdots \sigma'_1 (\sigma'_2\sigma'_1)^{k_2-1}(x) \rightarrow (\sigma'_2\sigma'_1)^{k_2}(x)=x.
 \end{align*}
 Note that $\sigma_1(a)=\sigma'_1(x)=a$, and that $\sigma_1$ and $\sigma'_1$, excluding $x$, differ only at the image of $a$. Similarly, $\sigma_2$ and $\sigma'_2$, excluding $x$, differ only at the image of $b$. Thus, the iterations starting with $\sigma_1(a)$ and $\sigma'_1(x)$ agree with each other until reaching $a$ or $b$.
 
 \emph{Claim~$1$.} The iteration starting with $x\rightarrow\sigma'_1(x)=a$ will not reach $a$ for a second time.\\
 This is because, otherwise, there must exists some $k\geq 1$, such that $\sigma'_1(\sigma'_2\sigma'_1)^k(x)=a$ or $(\sigma'_2\sigma'_1)^k(x)=a$.
 The former case can not happen, otherwise $(\sigma'_2\sigma'_1)^k(x)=x$, which will close the iteration instead of continuing to $a$.
 By Lemma~\ref{noodd}, the latter case, $(\sigma'_2\sigma'_1)^k(x)=\sigma'_2(\sigma'_1\sigma'_2)^{k-1}(\sigma'_1(x))=\sigma'_2(\sigma'_1\sigma'_2)^{k-1}(a)=a$ cannot happen either. Hence, Claim~$1$ follows.
 
 \emph{Claim~$2$.} The iteration starting with $\sigma'_1(x)=a$ will reach $b$ at least once.\\
 This is obvious since $\sigma'_1 (\sigma'_2\sigma'_1)^{k_2-1}(x)=b$.
 
 Now consider the first time the iteration starting with $\sigma'_1(x)$ reaches $b$.
 This must also be the first time the iteration starting with $\sigma_1(a)$ reaches $b$.
 There are two cases: either $(\sigma_2\sigma_1)^k(a)=b$ or $\sigma_1(\sigma_2\sigma_1)^k(a)=b$ for
 some $k\geq 1$.
 For the former case, we already have $a$ and $b$ as being in the same cycle of $\sigma_2\sigma_1$;\\
 For the latter case, we have $(\sigma_2\sigma_1)^{k+1}(a)=\sigma_2(b)=b$, which also implies that
 $a$ and $b$ are in the same cycle of $\sigma_2\sigma_1$.
 This completes the proof. \qed\\

 
 \begin{remark}
 	Lemma~\ref{keylem} can be alternatively proved in the following approach: first, we show that there is a way to assign signs `+' and `-' to elements in the set $T$ such that in both $\sigma_1$ and $\sigma_2$, every $2$-cycle has exactly one `+' element and one `-' elelment while $a$ and $b$ are positive, see a dihedral group action argument as in the Intersection-Theorem~\cite{reidys}. Then, we apply the Garsia-Milne Involution Principle~\cite{g-m-invo} to explain that $a$ and $b$ are in the same cycle of the product $\sigma_1 \sigma_2$.
 \end{remark}


\begin{lemma}\label{4lem1}
Let $\sigma_1,\: \sigma_2,\: \sigma'_1, \: \sigma'_2$ be defined as in the proof of Lemma~\ref{keylem}.
Then,
\begin{align}
C(\sigma'_1 \sigma'_2)-1= C(\sigma_1 \sigma_2).
\end{align}
\end{lemma}
\proof 
Following the discussion in the proof of Lemma~\ref{keylem}, 
any cycle not containing $a$ or $b$ or $x$ of $\sigma'_1 \sigma'_2$ is also
a cycle of $\sigma_1 \sigma_2$. Thus, the difference $C(\sigma'_1 \sigma'_2)-C(\sigma_1 \sigma_2)$
equals the difference of the number of cycles containing $a,\: b$ and $x$ in $\sigma'_1 \sigma'_2$
and the number of cycles containing $a$ and $b$ in $\sigma_1 \sigma_2$.

On the one hand, we have already shown that $a$ and $b$ are in the same cycle of $\sigma_1 \sigma_2$.
On the other hand, it is clear that $\sigma'_1 \sigma'_2(a)=b$; and by Claim~$1$ in the proof of Lemma~\ref{keylem}, $x$ and $a$
are not in the same cycle of $\sigma'_1 \sigma'_2$. Hence, the difference is exactly $1$, completing the proof. \qed

\begin{theorem} For any given signed permutation $a$,
\begin{align}
C(p \tilde{s})= C(\theta_1\theta_2)-1.
\end{align}
\end{theorem}
\proof 
By construction, the relation between the pair $\tilde{s}_{invo}, \: p_{invo}$ and $\theta_1, \: \theta_2$ is exactly the same as the pair $\sigma_1,\: \sigma_2$ and $\sigma'_1,\: \sigma'_2$.
Applying Lemma~\ref{4lem1}, we have
$$
C(\theta_1\theta_2)-1=C(p_{invo} \tilde{s}_{invo})=C(p \tilde{s}),
$$
 completing the proof. \qed

Accordingly, the lower bound on the reversal distance obtained by Chen and Reidys~\cite{chr-1} (and the topological-genus bound in the last section) is equivalent to the lower bound obtained by Bafna and Pevzner through breakpoint graphs~\cite{pev2,pev4}.

The discussion above has also implied Conjecture~\ref{7conj2}. Note that the conjecture can be equivalently formulated as follows:\\
Let $$p=(0,-1,-2,\ldots -n,n,n-1,\ldots 1)\:\text{ and }\:s=(0,s_1,s_2,\ldots s_n,-s_n,-s_{n-1},\ldots -s_1),$$
both being long cycles on the set $[n]^{\pm}$, and let $\pi=p\circ s$.
Then, the conjecture states that $n$ and $s_n$ are in the same cycle of $\pi$.

\begin{corollary}
	The elements $n$ and $s_n$ are in the same cycle of the product $p\circ s$.
\end{corollary}
\proof
Applying Lemma~\ref{keylem} to the involution decomposition 
$$p\circ s=p_{invo}\circ s_{invo}$$ 
the theorem follows.\qed

\end{document}